\documentclass[12pt]{article} 
\newcommand{\R}{\mathbf{R}}  
 \usepackage{latexsym}
   \font\twelvebm                       = cmmib10 at 12truept
   \font\tenbm                          = cmmib10 at 10truept
   \font\sevenbm                        = cmmib10 at 7truept
 = \tenbm \scriptscriptfont9              = \sevenbm
 
   at 10truept  at 10truept  at
 10truept  at 10truept
  at 10truept

 \mathchardef \BGamma            = "0900 \mathchardef \BDelta
 = "0901 \mathchardef \BTheta            = "0902 \mathchardef
 \BLambda           = "0903 \mathchardef \BXi               = "0904
 \mathchardef \BPi               = "0905 \mathchardef \BSigma
 = "0906 \mathchardef \BUpsilon          = "0907 \mathchardef \BPhi
 = "0908 \mathchardef \BPsi              = "0909 \mathchardef
 \BOmega            = "090A \mathchardef \Balpha            = "090B
 \mathchardef \Bbeta             = "090C \mathchardef \Bgamma
 = "090D \mathchardef \Bdelta            = "090E \mathchardef
 \Bepsilon          = "090F \mathchardef \Bzeta             = "0910
 \mathchardef \Beta              = "0911 \mathchardef \Btheta
 = "0912 \mathchardef \Biota             = "0913 \mathchardef
 \Bkappa            = "0914 \mathchardef \Blambda           = "0915
 \mathchardef \Bmu               = "0916 \mathchardef \Bnu
 = "0917 \mathchardef \Bxi               = "0918 \mathchardef \Bpi
 = "0919 \mathchardef \Brho              = "091A \mathchardef
 \Bsigma            = "091B \mathchardef \Btau              = "091C
 \mathchardef \Bupsilon          = "091D \mathchardef \Bphi
 = "091E \mathchardef \Bchi              = "091F \mathchardef \Bpsi
 = "0920 \mathchardef \Bomega            = "0921 \mathchardef
 \Bvarepsilon       = "0922 \mathchardef \Bvartheta         = "0923
 \mathchardef \Bvarpi            = "0924 \mathchardef \Bvarrho
 = "0925 \mathchardef \Bvarsigma         = "0926 \mathchardef
 \Bvarphi           = "0927
 \mathchardef \bA        = "0941 \mathchardef \bB        = "0942
 \mathchardef \bC        = "0943 \mathchardef \bD        = "0944
 \mathchardef \bE        = "0945 \mathchardef \bF        = "0946
 \mathchardef \bG        = "0947 \mathchardef \bH        = "0948
 \mathchardef \bI        = "0949 \mathchardef \bJ        = "094A
 \mathchardef \bK        = "094B \mathchardef \bL        = "094C
 \mathchardef \bM        = "094D \mathchardef \bN        = "094E
 \mathchardef \bO        = "094F \mathchardef \bP        = "0950
 \mathchardef \bQ        = "0951 \mathchardef \bR        = "0952
 \mathchardef \bS        = "0953 \mathchardef \bT        = "0954
 \mathchardef \bU        = "0955 \mathchardef \bV        = "0956
 \mathchardef \bW        = "0957 \mathchardef \bX        = "0958
 \mathchardef \bY        = "0959 \mathchardef \bZ        = "095A
 \mathchardef \ba        = "0961 \mathchardef \bb        = "0962
 \mathchardef \bc        = "0963 \mathchardef \bd        = "0964
 \mathchardef \bee       = "0965 
 \mathchardef \bff       = "0966 \mathchardef \bg        = "0967
 \mathchardef \bh        = "0968
 \mathchardef \bj        = "096A \mathchardef \bk        = "096B
 \mathchardef \bl        = "096C \mathchardef \bm        = "096D
 \mathchardef \bn        = "096E \mathchardef \bo        = "096F
 \mathchardef \bp        = "0970 \mathchardef \bq        = "0971
 \mathchardef \br        = "0972 \mathchardef \bs        = "0973
 \mathchardef \bt        = "0974 \mathchardef \bu        = "0975
 \mathchardef \bv        = "0976 \mathchardef \bw        = "0977
 \mathchardef \bx        = "0978 \mathchardef \by        = "0979
 \mathchardef \bz        = "097A

 \font\tencb            = cmssbx10 scaled \magstep4 \font\eigcb
 = cmssbx10 scaled \magstep2 \textfont8             = \tencb
 \mathchardef\bAs       = "1841
 \def\Asem#1#2{\mathop{\vrule height10.5pt depth5.5pt width0pt\bAs}_{#1}^{#2}}
 \def\asem#1#2{
          \ifmmode
         \ifinner
            \raise0.9pt\hbox{$\scriptstyle\bAs$}_{#1}^{#2}
         \else
            \Asem{#1}{#2}
         \fi
          \fi
          }

 \newtheorem{theo}{\small\bf Theorem}
 \newtheorem{lem}{\small\bf Lemma}
 \newtheorem{prop}{\small\bf Proposition}
 \newtheorem{rem}{\small\bf Remark}
 \newenvironment{REM}{\begin{rem} \rm}{\end{rem}}
 \newtheorem{exam}{\small\bf Example}
 \newenvironment{EXAM}{\begin{exam} \rm}{\end{exam}}
  \newtheorem{defi}{\small\bf Definition}
 \newenvironment{DEFI}{\begin{defi} \rm}{\end{defi}}
 \renewcommand{\Pr}{\mbox{\rm  \hspace*{.2ex}I\hspace{-.5ex}P\hspace*{.2ex}}}

 \newcommand{\be}{\begin{equation}}
 \newcommand{\ee}{\end{equation}}

 \newcommand{\E}{\mbox{\rm \hspace*{.2ex}I\hspace{-.5ex}E\hspace*{.2ex}}}

 \newcommand{\Var}{\mbox{\rm \hspace*{.2ex}Var\hspace*{.2ex}}}

 \newcommand{\bbb}[1]{\mbox{\boldmath $ #1 $}}

 \newenvironment{pr}[1]{{\small\bf {#1}:}}{}

 \title{ \Large\bf Finite sequences
 representing expected order statistics}
 \author{\normalsize\bf
 {\rm by}
 \\
 A.\  Okolewski and N.\ Papadatos
 }
 \date{\small
 Institute of Mathematics,
 Lodz University of Technology,
 90-924 Lodz, Poland
 \\
 and
 \\
 Department of Mathematics,
 National and Kapodistrian
 University of Athens,  Panepistemiopolis, 157 84
 Athens,
 Greece.}

\begin{document}

 \maketitle
 \vspace*{-2em}

 \begin{abstract}
 Characterizations of finite
 sequences
 $\beta_{1}<\cdots<\beta_{n}$
 representing expected values of order statistics from a random sample
 of size $n$ are given. As a by-product, a characterization of binomial
 mixtures, when the mixing random variable is supported in the
 open interval $(0,1)$, is presented; this enables the exact
 description of the convex hull of the open binomial curve, as
 well as the open moment curve.
 \vspace{3ex}


 \end{abstract}

 \thispagestyle{empty}

 \section{Introduction}
 \label{sec.intro}
 In the present note we consider the following problem:
 Given $n$ real numbers
 \[
 \beta_{1}<\cdots<\beta_{n},
 \]
 under what conditions on $\beta$'s is there an integrable
 random variable (r.v.)\ $X$
 such that
 \[
 \E X_{i:n}=\beta_{i}, \ \ 1\leq i\leq n ?
 \]
 [Here, $X_{1:n}\leq \cdots\leq X_{n:n}$ are the order statistics of
 independent, identically distributed r.v.'s $X_1,\ldots,X_n$,
 each with distribution like $X$.] \
 Notice that
 the number $n$ is held fixed; the
 question for infinite sequences is closely connected to the Hausdorff (1921)
 moment problem, and its answer is well-known
 from the works of  Huang (1998), Kadane (1971, 1974), Kolodynski (2000) and
 Papadatos (2017). Some relative results for the finite case can be found
 in Mallows (1973).

 \section{A characterization of finite sequences of expected order statistics}
 \label{sec.2}
 Without loss of generality we may consider the numbers
 \[
 \widetilde{\beta}_{i} =\frac{\beta_{i}-c}{\lambda}, \ \
 c\in\R, \ \lambda>0,
 \]
 instead of $\beta_{i}$. Clearly these numbers will be the expected
 order statistics
 (=EOS) from $(X-c)/\lambda$ if and only if the $\beta$'s are the EOS from
 $X$.

 \noindent
 First, we seek for a necessary condition. Assume that $X$
 is a non-degenerate
 random variable
 with distribution
 function (d.f.) $F$ and
 $\E|X|<\infty$. Let $X_1,\ldots,X_n$
 be independent, identically distributed (i.i.d.)
 random variables
 with d.f.\ $F$, and denote by $X_{1:n}\leq \cdots\leq X_{n:n}$
 the corresponding order statistics.
 It is known that
 \be
 \label{muk from mjn}
  \E X_{k:k}:=\mu_k=
   {{n\choose k}}^{-1}\sum_{j=k}^n {j-1 \choose k-1} \mu_{j:n},
   \ \ \ k=1,\ldots,n,
 \ee
 where $\mu_{j:n}:=\E X_{j:n}$, $j=1,\ldots,n$;
 this follows
 by a trivial application of Newton's formula
 to the expression
 $
 \mu_k
 =
  k \int_{0}^1 u^{k-1} F^{-1}(u) [u+(1-u)]^{n-k} du,
 $
 where
 $F^{-1}(u):=\inf\{x:F(x)\geq u\}$, $0<u<1$,
 is the left-continuous inverse of $F$.
 From (\ref{muk from mjn}) with $k=1,2$,
 \be
 \label{mu1 mu2}
 \mu_1=\frac{\mu_{1:n}+\cdots+\mu_{n:n}}{n}, \ \ \
 \mu_2=\frac{2}{n(n-1)}\sum_{j=2}^n (j-1)\mu_{j:n}.
 \ee
 On the other hand,
 it is well-known (see Jones and Balakrishnan (2002))
 that
 \be
 \label{Pearson}
 \mu_{j+1:n}-\mu_{j:n} ={n\choose j} \int_{\alpha}^\omega F(x)^j
 (1-F(x))^{n-j} dx>0, \ \ \ j=1,\ldots, n-1,
 \ee
 where $\alpha<\omega$ are the
 endpoints
 of the support of
 $X$;
 actually this formula goes back to Karl Pearson (1902).
 Notice that $-\infty\leq \alpha<\omega\leq \infty$,
 $\alpha<\omega$ because $F$ is non-degenerate, and the
 integral in (\ref{Pearson}) is finite since $X$ is
 integrable.
 From (\ref{muk from mjn}) and
 (\ref{Pearson}) (applied to $n=2$),
 \[
 \mu_2-\mu_1=\int_{\alpha}^\omega F(x)(1-F(x)) dx,
 \]
 while (\ref{mu1 mu2}) yields
 \be
 \label{mu2-mu1}
 \mu_2-\mu_1
 =\frac{1}{n(n-1)}\sum_{i=1}^{n-1} i(n-i)(\mu_{i+1:n}-\mu_{i:n}).
 \ee
 Choosing $c=n^{-1}\sum_{j=1}^n \mu_{j:n}$,
 $\lambda
 =\big(n(n-1)\big)^{-1}\sum_{i=1}^{n-1} i(n-i)(\mu_{i+1:n}-\mu_{i:n})>0$,
 the numbers $\widetilde{\mu}_{j:n}=(\mu_{j:n}-c)/\lambda$ are the
 EOS from $(X-c)/\lambda$. Therefore, by considering
 $(X-c)/\lambda$ in place of $X$,
 we may further assume that $\mu_2-\mu_1=1$. Then,
 \[
 \int_{\alpha}^\omega F(x) (1-F(x)) dx =1.
 \]
 Since $F(y)(1-F(y))>0$ for $y\in (\alpha,\omega)$, and zero outside
 $[\alpha,\omega)$, it follows that
 $f_Y(y):=F(y)(1-F(y))$ defines a Lebesgue density of a random variable, say $Y$,
 supported in the (finite or infinite) interval $(\alpha,\omega)$.
 From (\ref{mu2-mu1}),
 the numbers $p_j:= j(n-j)(\mu_{j+1:n}-\mu_{j:n})/(n(n-1))$,
 $j=1,\ldots,n-1$,
 are strictly positive probabilities (summing to $1$).
 Also, the integral
 in (\ref{Pearson}) can be rewritten as
 \[
 \int_{\alpha}^\omega \Big(F(y)^{j-1}
 (1-F(y))^{n-j-1}\Big) f_Y(y) dy=\E\Big\{ T^{j-1} (1-T)^{n-j-1}\Big\},
 \]
 where $T:=F(Y)$ is a random variable taking values in the
 interval $(0,1)$ w.p.\ 1, because, by definition,
 $\Pr(\alpha<Y<\omega)=1$. Hence,
 (\ref{Pearson}) reads as
 \[
 p_j=\frac{j(n-j)}{n(n-1)}{n\choose j}\E\Big\{ T^{j-1} (1-T)^{n-j-1}\Big\},
 \ \ \ j=1,\ldots,n-1,
 \]
 and we have shown the following
 \begin{prop}
 \label{prop.1}
 If $X_1,\ldots,X_n$ are i.i.d.\ integrable non-degenerate r.v.'s, then
 there exists an r.v.\ $T$, with $\Pr(0<T<1)=1$, such that
 \begin{eqnarray}
  \nonumber
  \frac{(j+1)(n-j-1)(\mu_{j+2:n}-\mu_{j+1:n})}
  {\sum_{i=1}^{n-1}i(n-i)(\mu_{i+1:n}-\mu_{i:n})}
  \hspace*{-1.4ex}
  &
  =
  &
  \hspace*{-1.4ex}
  \E\left\{ {n-2\choose j} T^{j} (1-T)^{n-2-j}\right\},
  \\
  &&
  \label{binomial moments}
 \ \ \ \ \ \ \ \ \ \ \ \ \ \ \ \ j=0,\ldots,n-2.
 \end{eqnarray}
 \end{prop}

 \noindent
 It is of interest to observe that the
 binomial moments of $T$ appear in the r.h.s.\ of (\ref{binomial moments}).
 Clearly, the r.v.\ $T$ in this representation need not be unique; any other
 r.v.\ $T'$ with
 $\Pr(0<T'<1)=1$, possessing
 identical moments up to order $n-2$ with $T$, will fulfill the
 same relationship.

 \begin{REM}
 \label{rem.1}
 For any integrable non-degenerate r.v.\
 $X$ with d.f.\ $F$ we may define the r.v.\ $T$ as in the proof
 of Proposition \ref{prop.1}, that is, $T=F(Y)$ where $Y$ has density
 $f_Y(y)=F(y)(1-F(y))/\lambda$ with $\lambda=\int
 F (1-F)$.
 It can be shown,
 using Lemma 4.1 in Papadatos (2001),
 that the d.f.\ of $T$ is specified by
 \be
 \label{FT}
 \Pr(T<t)=\frac{1}{\lambda}
 \Bigg[t(1-t)F^{-1}(t)-\int_{0}^t (1-2u) F^{-1}(u) d u\Bigg],
 \ \ 0<t<1.
 \ee
 Notice that $\lambda=\int_{0}^1 (2t-1)F^{-1}(t) dt$ and, hence,
 the function $F^{-1}$ determines uniquely the d.f.\ of $T$.
 Moreover, (\ref{FT}) shows that the entire location-scale family of $X$,
 $\{c+\lambda X: \ c\in\R, \ \lambda>0\}$,
 is mapped to a single r.v.\ $T\in(0,1)$.
 Provided
 that $X$ has (finite or infinite) interval support, non-vanishing
 density $f$ and differentiable inverse d.f.\ $F^{-1}$, we
 conclude from (\ref{FT}) that a density of $T$ is given by
 \be
 \label{density fT}
 f_T(t)=\frac{t(1-t)}{\lambda f(F^{-1}(t))}, \ \ 0<t<1.
 \ee
 \end{REM}

 Next, we proceed to verify that the preceding procedure can be inverted,
 showing sufficiency of (\ref{binomial moments}).
 To this end, we shall
 make use of the following lemma, which is of independent
 interest in itself. A detailed proof is postponed to the appendix.

 \begin{lem}
 \label{lem.1}
 Let $T$ be an r.v.\ with d.f.\ $F_T$ such that $\Pr(0<T<1)=1$.
 Then, there exists a unique, non-degenerate, integrable,
 r.v.\ $X$,
 satisfying
 \be
 \label{EXkk}
 \E X_1 =0  \ \ \mbox{and} \ \ \ \E X_{k+2:k+2}-\E X_{k+1:k+1}=\E T^k,
 \ \ k=0,1,\ldots,
 \ee
 where $X_{k:k}=\max\{X_1,\ldots,X_k\}$ with
 $X_1,X_2,\ldots$ being i.i.d.\ copies of $X$.
 The inverse distribution function of $X$ is given
 by
 \be
 \label{DistributionInverse}
 F_0^{-1}(t) :=\frac{F_T(t-)}{t(1-t)}
 - 4 F_T(\mbox{$\frac{1}{2}-$})
 -
 \int_{1/2}^{t} \frac{2u-1}{u^2(1-u)^2} F_T(u)d u - c_T,
 \ee
 $0<t<1$,
 where $F_T(t-)=\Pr(T<t)$, $\int_{1/2}^t du=-\int_{t}^{1/2} du$ for $t<1/2$,
 \be
 \label{cT}
 c_T := \E\Bigg[\frac{1}{T}I(T\geq \mbox{$\frac{1}{2}$})\Bigg]-
 \E\Bigg[\frac{1}{1-T}I(T<\mbox{$\frac{1}{2}$})\Bigg],
 \ee
 and $I$ denotes an indicator function.
 \end{lem}

 \begin{REM}
 \label{rem.2}
 Any r.v.\ $T\in(0,1)$ can be viewed as the {\it expected order statistics
 generator} of its corresponding r.v.\ $X$ with inverse d.f.\ $F_0^{-1}$
 as in (\ref{DistributionInverse}). This is so because
 the map $T\rightarrow X$ (i.e., $F_T\rightarrow F_0\equiv F_X$),
 defined implicitly by Lemma \ref{lem.1},
 is
 one to one and onto
 from the space ${\cal T}=\{T: \Pr(0<T<1)=1\}$
 to
 ${\cal H}=\{X: \E X=0, \E X_{2:2}=1\}$, where
 identically distributed r.v.'s are considered as equal.
 Its inverse is given
 by Remark \ref{rem.1} (with $\lambda=1$, since $X\in {\cal H}$).
 In view of (\ref{mujn from muk}), below, it is the suitable
 (and unique) transformation
 that quantifies the
 characterization of Hoeffding (1953), stating that the sequence of
 expected order statistics
 characterizes the corresponding distribution.
 It also provides an explicit connection of the (infinite) sequence of
 expected order statistics to the Hausdorff (1921) moment problem;
 see Kadane (1971, 1974), Huang (1998), Kolodynski (2000), Papadatos
 (2017).
 \end{REM}

 \begin{REM}
 \label{rem.3}
 Suppose that the r.v.\ $T$ of Lemma \ref{lem.1} is absolutely
 continuous with density $f_T$.
 Assume also that the corresponding r.v.\ $X$
 (with $\E X=0$, $\E X_{2:2}=1$, inverse
 d.f.\
 $F_0^{-1}$ as in (\ref{DistributionInverse}))
 is absolutely continuous, admitting a non-vanishing density $f_0$ in the
 (finite or infinite)
 interval support of $X$, and that $F_0^{-1}$ is differentiable. Then
 (see Remark \ref{rem.1}),
 \[
 f_T(t)=\frac{t(1-t)}{f_0(F_0^{-1}(t))}, \ \ \
 F_0^{-1}(t)=\int_{1/2}^t\frac{f_T(u)}{u(1-u)} du -c_T,
 \ \ \ 0<t<1.
 \]
 For example, if $T$ is Beta$(2,2)$ then $X$ is uniform in $(-3,3)$;
 if $T$ is Beta$(2,1)$ then $X=2{\cal E}-2$ where ${\cal E}$ is standard exponential;
 if $T$ is Beta$(1,2)$ then $X=2-2{\cal E}$; if
 $T$ is standard uniform then $X$ is standard logistic with density
 $f_0(x)=e^{-x}/(1+e^{-x})^2$, $x\in\R$; if $T$ is degenerate with
 $\Pr(T=\rho)=1$ then (\ref{DistributionInverse}) shows that $X$ is a two-valued
 r.v.\
 with $\Pr(X=-1/\rho)=\rho$, $\Pr(X=1/(1-\rho))=1-\rho$.
 \end{REM}

 The characterization for finite $n$ reads as follows.

 \begin{theo}
 \label{theo.1}
 Given $n$ real numbers $\beta_{1}<\cdots<\beta_{n}$,
 the following
 are equivalent.
 \smallskip

 \noindent
 {\rm (i)}
 The $\beta$'s are EOS, that is, there exist i.i.d.\ integrable non-degenerate
 r.v.'s \\
 $X_1,\ldots,X_n$
 such that $\E X_{j:n}=\beta_{j}$, $j=1,\ldots,n$.
 \smallskip

 \noindent
 {\rm (ii)}
 There exists an r.v.\ $T$, with $\Pr(0<T<1)=1$,
 such that
 \begin{eqnarray}
 \nonumber
   \frac{(j+1)(n-j-1)(\beta_{j+2}-\beta_{j+1})}
   {\sum_{i=1}^{n-1}i(n-i)(\beta_{i+1}-\beta_{i})}
  \hspace*{-1ex}&=&\hspace*{-1ex}
  \E\left\{ {n-2\choose j} T^{j} (1-T)^{n-2-j}\right\},
  \\
  \label{betaj binomial}
  &&
 \ \ \ \ \ \ \ \ \ \ \ \ \ \ j=0,\ldots,n-2.
 \end{eqnarray}

 \noindent
 {\rm (iii)}
 There exists an r.v.\ $T$, with $\Pr(0<T<1)=1$,
 such that
 \begin{eqnarray}
 \nonumber
 \frac{n-1}{{n-1 \choose k+1}\sum_{i=1}^{n-1}i(n-i)(\beta_{i+1}-\beta_i)}
 \sum_{j=k+1}^{n-1}(n-j){j\choose k+1}(\beta_{j+1}-\beta_j)
  \hspace*{-1ex}&=&\hspace*{-1ex}
  \E T^k,
  \\
  \label{simple.moments}
 k=0,\ldots,n-2.
 \end{eqnarray}
 \end{theo}
 \begin{pr}{Proof} The equivalence of (\ref{betaj binomial}) and (\ref{simple.moments})
 follows by a straightforward computation, while the implication
 (i)$\Rightarrow$(ii) is proved in Proposition \ref{prop.1}.
 In order to verify (ii)$\Rightarrow$(i), assume
 that (\ref{betaj binomial}) is satisfied for some $T$ with $\Pr(0<T<1)=1$,
 and consider the r.v.\ $X$ as defined in Lemma \ref{lem.1}.
 Let $\mu_{j:n}=\E X_{j:n}$ and $\mu_k=\E X_{k:k}$. Then,
 \be
 \label{mujn from muk}
 \mu_{j:n}=n {n-1 \choose j-1} \sum_{i=j}^n (-1)^{i-j} {n-j \choose i-j} \frac{\mu_i}{i},
 \ \ \ j=1,\ldots,n;
 \ee
 see Mallows (1973), Arnold {\it et al.}\ (1992), David and Nagaraja (2003).
 It follows that
 \[
 \mu_{j+2:n}-\mu_{j+1:n}={n \choose j+1} \sum_{i=j+1}^n
 (-1)^{i-j} {n-j-1 \choose i-j-1} \mu_i,
 \ \ \ j=0,\ldots,n-2.
 \]
 By a trivial application of the binomial theorem to $(1-T)^{n-2-j}$, and since
 $\E T^k=\mu_{k+2}-\mu_{k+1}$, see (\ref{EXkk}), we obtain
 \begin{eqnarray*}
 &&  \E\left\{ {n-2\choose j} T^{j} (1-T)^{n-2-j}\right\}
  = {n-2\choose j} \sum_{i=j+1}^n (-1)^{i-j} {n-j-1 \choose i-j-1} \mu_i,
  \\
  && \hspace*{50ex} j=0,\ldots,n-2.
 \end{eqnarray*}
 Hence, for $j=0,\ldots,n-2$,
 \[
 \frac{(j+1)(n-j-1)}{n(n-1)}\Big(\mu_{j+2:n}-\mu_{j+1:n}\Big)=
 \E\left\{ {n-2\choose j} T^{j} (1-T)^{n-2-j}\right\},
 \]
 and (\ref{betaj binomial}) implies  that  for some $\lambda>0$,
 \[
 \beta_{j+1}-\beta_{j}=\lambda (\mu_{j+1:n}-\mu_{j:n}),
 \ \ j=1,\ldots, n-1.
 \]
 It follows by induction on $j$ that $\beta_{j}=\beta_{1}+\lambda(\mu_{j:n}-\mu_{1:n})$,
 and therefore, $(\beta_{j}-c)/\lambda=\mu_{j:n}$, with $c=\beta_{1}-\lambda \mu_{1:n}$.
 Hence, the
 numbers $\big((\beta_{j}-c)/\lambda\big)_{j=1}^n$ are expected order statistics,
 and thus, the same is true for $\beta$'s. 
 \end{pr}

 \begin{REM}
 \label{rem.new}
 The r.h.s.\ of (\ref{betaj binomial}) corresponds to a
 {\it Binomial Mixture} (of a particular form, since $\Pr(T=0)=\Pr(T=1)=0$).
 Obviously, the necessary and sufficient condition
 (\ref{simple.moments}) is always satisfied for $n=2$ and $n=3$.
 Hence, the true
 problem begins at $n=4$.
 \end{REM}

 \section{The truncated moment problem for finite open intervals}
 \label{sec.3}
 In this section, we obtain a precise characterization by invoking
 results from the {\it truncated moment problem for finite intervals}.
 The
 existing results are limited to
 compact intervals
  and are not
  applicable
 to our case, since, according to the characterization of
 Theorem \ref{theo.1},
 a suitable $T$ lies in the {\it open interval} $(0,1)$ w.p.\ 1.

 \begin{DEFI}
 \label{defi.nu_k}
 Given $n\geq 4$ numbers $\beta_1<\cdots<\beta_n$,
 let ${\bbb \beta}=(\beta_1,\ldots,\beta_n)$
 and define the vector $(\nu_k)_{k=0}^{n-2}={\bbb \nu}={\bbb \nu}({\bbb \beta})$
 by
 \be
 \label{nu_k}
 \nu_k:=\frac{n-1}{\lambda {n-1 \choose k+1}}
 \sum_{j=k+1}^{n-1}(n-j){j\choose k+1}(\beta_{j+1}-\beta_j), \ \ \ \ k=0,\ldots,n-2,
 \ee
 where
 $\lambda=\lambda(\bbb \beta):=\sum_{i=1}^{n-1}i(n-i)(\beta_{i+1}-\beta_i)>0$.
 \end{DEFI}

 \noindent
 It is easily checked that $1=\nu_0>\nu_1>\cdots>\nu_{n-2}>0$, and that
 the vector $\bbb \nu$ is invariant under location-scale transformations on
 the $\beta$'s.

 According to Theorem \ref{theo.1}, the $\beta$'s are EOS if and only if
 the $\nu$'s fulfill the truncated moment problem in the interval
 $(0,1)$. However, for the truncated moment problem,
 well-known results exist for a compact interval
 $[a,b]$; see, e.g., Theorem IV.1.1 of Karlin and Studden (1966) or
 Theorems 10.1, 10.2 in Schm\"{u}dgen (2017).
 In order to obtain the corresponding necessary and sufficient conditions
 for open intervals, we shall make use
 of the following
 \begin{theo}
 \label{theo.Sch}
 {\rm (Richter-Tchakaloff Theorem; see Schm\"{u}dgen (2017), Theorem 1.24).}
 Let $({\mathcal X},{\mathcal F},\mu)$ be a measure space
 and $V$ a finite dimensional linear subspace of $L^1(\R,\mu)$.
 Define the linear functional $L_{\mu}$ by $L_{\mu}(f):=\int f d\mu$, $f\in V$.
 Then, there exists a measure $\mu_0$ in $({\mathcal X},{\mathcal F})$,
 supported on $k\leq \dim V$ points
 of $\mathcal X$, such that $L_{\mu_0}\equiv L_{\mu}$ on $V$, that is,
 $\int f d\mu_0=\int f d\mu$ for all $f\in V$.
 \end{theo}

 \noindent
 A symmetric $n\times n$ matrix $A$ with real entries is positive definite
 (denoted by $A \succ 0$) if ${\bbb x}^T A{\bbb x}>0$ for all
 ${\bbb x}\in \R^n\setminus\{{\bbb 0}\}$, where
 ${\bbb x}^T$ denotes the transpose of a column vector ${\bbb x}\in \R^n$.
 Similarly, $A$ is positive semi-definite (or nonnegative definite) if
 ${\bbb x}^T A{\bbb x}\geq 0$ for all ${\bbb x}\in \R^n$, and this is denoted
 by $A\succeq 0$.

 \begin{DEFI}
 \label{defi.Hankel}
 (Hankel matrices).
 Let $n\in\{4,5,\ldots\}$, $0\leq\varepsilon< 1/2$,
 and consider the numbers $\nu_k$ as in Definition
 \ref{defi.nu_k}.
 \smallskip

 \noindent
 (i) Case $n=2m+2$: We define
 \be
 \label{Hankel.even}
 A_0(\varepsilon):=\Big(\nu_{i+j}\Big)_{i,j=0}^m, \ \ \
 B_0(\varepsilon):=\Big(\nu_{i+j+1}-\nu_{i+j+2}
 -\varepsilon(1-\varepsilon) \nu_{i+j}\Big)_{i,j=0}^{m-1},
 \ee
 and $A_0:=A_0(0)$, $B_0:=B_0(0)$.
 \smallskip

 \noindent
 (ii) Case $n=2m+3$: We define
 \be
 \label{Hankel.odd}
 A_1(\varepsilon):=\Big( \nu_{i+j+1}-\varepsilon\nu_{i+j}\Big)_{i,j=0}^m, \ \ \
 B_1(\varepsilon):=\Big((1-\varepsilon)\nu_{i+j}-\nu_{i+j+1}\Big)_{i,j=0}^{m},
 \ee
 and $A_1:=A_1(0)$, $B_1:=B_1(0)$.
 \end{DEFI}

 Notice that the matrices $A_0(\varepsilon),
 A_1(\varepsilon), B_1(\varepsilon)$ are of
 order $m+1$, while
 $B_0(\varepsilon)$ is of order $m$.
 The following theorem contains our main result;
 compare with Mallows (1973).

 \begin{theo}
 \label{theo.2}
 Let $n\in \{4,5,\ldots\}$,
 $\beta_1<\cdots<\beta_n$, and $(\nu_0,\ldots,\nu_{n-2})$
 as in Definition
 \smallskip
 {\rm \ref{defi.nu_k}}.

 \noindent
 {\rm{(i)}} If $n=2m+2,$ then
 the $\beta$'s are EOS if and only if
 $A_0(\varepsilon)\succeq 0$ and $B_0(\varepsilon)\succeq 0$
 for some $\varepsilon\in(0,1/2)$, where $A_0(\varepsilon)$ and $B_0(\varepsilon)$
 are given by Definition
 \smallskip
 {\rm \ref{defi.Hankel}(i)}.

 \noindent
 {\rm{(ii)}} If $n=2m+3,$
 then
 the $\beta$'s are EOS if and only if
 $A_1(\varepsilon)\succeq 0$ and $B_1(\varepsilon)\succeq 0$
 for some $\varepsilon\in(0,1/2)$, where $A_1(\varepsilon)$ and $B_1(\varepsilon)$
 are given by Definition {\rm \ref{defi.Hankel}(ii)}.
 \smallskip

 \noindent
 {\rm{(iii)}} If $n=2m+2$, the
 condition $A_0\succ 0$ and $B_0\succ 0$ is sufficient, but not necessary, for the
 $\beta$'s to be EOS. Similarly,
 if $n=2m+3$, the
 condition $A_1\succ 0$ and $B_1\succ 0$ is sufficient, but not necessary, for the
 $\beta$'s to be EOS.
 \end{theo}
 \begin{pr}{Proof}
 (i) and (ii): According to Theorems 10.1, 10.2
 in Schm\"udgen (2017),
 or
 Theorem IV.1.1 of Karlin and Studden (1966), the
 condition $A_i(\varepsilon)\succeq 0$ and $B_i(\varepsilon)\succeq 0$
 ($i=0$ or $1$) is necessary and sufficient for $(\nu_k)_{k=0}^{n-2}$ to be
 a truncated moment sequence in the interval $[\varepsilon,1-\varepsilon]$.

 Assume first that $A_i(\varepsilon)\succeq 0$ and $B_i(\varepsilon)\succeq 0$
 for some
 $\varepsilon\in(0,1/2)$.
 Since
 $\nu_0=1$, any solution (=representing measure) $\mu$ will
 be a probability measure. Equivalently,
 the r.v.\ $T$ with d.f.\ $F_T(x)=\mu\big((-\infty,x]\big)$ takes values
 in $[\varepsilon,1-\varepsilon]\subseteq(0,1)$ and satisfies
 $\E T^k=\nu_k$, $k=0,\ldots,n-2$. From Theorem
 \ref{theo.1}(iii) it follows that the $\beta$'s are EOS.

 To prove necessity, assume that the $\beta$'s are EOS. From
 Theorem  \ref{theo.1}(iii) we can find an r.v.\ $T$ with
 $\E T^k=\nu_k$, $k=0,\ldots,n-2$, and $\Pr(0<T<1)=1$.
 Let $\mu_T$ be the probability measure of $T$ and consider
 the probability space $({\mathcal X},{\mathcal F},\mu):=
 ((0,1),{\mathcal B},\mu_T)$, where ${\mathcal B}$
 is the Borel $\sigma-$field on $(0,1)$.
 Define the
 space $V$ of real polynomials $f:(0,1)\rightarrow\R$ of degree
 $\leq n-2$;
 obviously, $V$ is a linear subspace of $L^1(\R,\mu_T)$
 of dimension $n-1$ (finite). Consider also the
 Riesz functional $L_{\mu_T}:V\rightarrow\R$
 defined by $L_{\mu_T}(f):=\int f d \mu_T=\sum_{k=0}^{n-2}a_k \nu_k$
 for $f(x)=\sum_{k=0}^{n-2}a_k x^k \in V$.
 Form Richter-Tchakaloff Theorem (see Theorem \ref{theo.Sch}, above),
 there exists a measure $\mu_0$, supported in at most $n-1$ points
 of ${\mathcal X}=(0,1)$, such that $L_{\mu_0}\equiv L_{\mu_T}$ on $V$;
 in particular, $\nu_k=\int_{(0,1)}x^k d {\mu_T}(x)=\int_{(0,1)}x^k d {\mu_0}(x)$,
 $k=0,\ldots,n-2$. Thus, $\mu_0$ is a probability measure ($\nu_0=1$)
 supported on a finite number of points in $(0,1)$, possessing  the same
 initial $n-2$ moments as $\mu_T$. This means that $\mu_0$ solves the truncated
 moment problem for $(\nu_k)_{k=0}^{n-2}$
 in the interval $[t_1,t_{2}]$, where $t_1\in(0,1)$ is the minimum
 supporting point of $\mu_0$ and $t_{2}\in(0,1)$ the maximum one.
 Choose $\varepsilon>0$
 such that $\varepsilon<\min\{t_1,1-t_{2}\}$. Then,
 the sequence $(\nu_k)_{k=0}^{n-2}$ is the moment sequence
 of $\mu_0$, supported in the interval $[\varepsilon,1-\varepsilon]$,
 and  Theorems 10.1, 10.2 in Schm\"udgen (2017)
 imply that $A_i(\varepsilon)\succeq 0$ and $B_i(\varepsilon)\succeq 0$
 ($i=0$ or $1$).
 \smallskip

 \noindent {\rm (iii)}
 First we prove sufficiency.
 Denote by $\lambda_{\min}(M)$ (resp.\ $\lambda_{\max}(M)$) the smallest
 (resp.\ the largest) eigenvalue of a real symmetric matrix $M$.
 For the case $n=2m+2$, the matrix $A_0(\varepsilon)$
 is independent of $\varepsilon$, hence, $A_0(\varepsilon)=A_0\succ 0$ by hypothesis.
 Moreover,
 $B_0(\varepsilon)=
 B_0-\varepsilon(1-\varepsilon)M_0$ for some real symmetric matrix $M_0$;
 see (\ref{Hankel.even}).
 Since $\lambda_{\min}(B_0)>0$ by assumption,
 it follows
 that for any ${\bbb x}=(x_0,\ldots,x_{m-1})^T\in\R^{m}$,
 ${\bbb x}^T B_0(\varepsilon) {\bbb x}={\bbb x}^T B_0 {\bbb x}
 -\varepsilon(1-\varepsilon) {\bbb x}^T M_0 {\bbb x}\geq
 \big[\lambda_{\min}(B_0)-\varepsilon(1-\varepsilon)\lambda_{\max}(M_0)\big]{\bbb x}^T
 {\bbb x}\geq 0$, if $\varepsilon >0$ is sufficiently small.
 Hence, the sufficient condition (i), namely, $A_0(\varepsilon)\succeq 0$
 and $B_0(\varepsilon)\succeq 0$ for some small $\varepsilon>0$, is satisfied.
 Similarly, when $n=2m+3$ we have
 $A_1(\varepsilon)=A_1-\varepsilon M_1$ and
 $B_1(\varepsilon)=B_1-\varepsilon M_1$
 for some real symmetric matrix  $M_1$;
 see (\ref{Hankel.odd}).
 From $\lambda_{\min}(A_1)>0$, $\lambda_{\min}(B_1)>0$,
 it follows that for any ${\bbb x}\in \R^{m+1}$,
 ${\bbb x}^T A_1(\varepsilon){\bbb x}={\bbb x}^T A_1 {\bbb x}
 -\varepsilon {\bbb x}^T M_1 {\bbb x}\geq \big[\lambda_{\min}(A_1)
 -\varepsilon\lambda_{\max}(M_1)\big]{\bbb x}^T
 {\bbb x}\geq 0$, and
 ${\bbb x}^T B_1(\varepsilon){\bbb x}\geq
 \big[\lambda_{\min}(B_1)-\varepsilon\lambda_{\max}(M_1)\big]{\bbb x}^T
 {\bbb x}\geq 0$, provided $\varepsilon>0$ is sufficiently small.
 Hence, the sufficient condition (ii),
 $A_1(\varepsilon)\succeq 0$
 and $B_1(\varepsilon)\succeq 0$ for some small $\varepsilon>0$,
 is satisfied. Therefore,
 in both cases, the condition (iii) is sufficient for the
 $\beta$'s to represent EOS.

 Finally, we show that the condition (iii),
 namely $A_i\succ 0$ and $B_i\succ 0$ ($i=0$ or $1$),
 is not necessary. To this end, consider the
 sequence
 $\beta_j:=\sum_{k=n+1-j}^n {n \choose k}$, $j=1,\ldots,n$.
 Then, $\beta_{j+1}-\beta_j={n\choose j}$ ($j=1,\ldots,n-1$)
 and a straightforward computation
 yields $\nu_k=2^{-k}$, $k=0,\ldots,n-2$; see (\ref{nu_k}).
 Suppose first that $n=2m+2$ and let ${\bbb x}^T=(x_0,\ldots,x_m)\in \R^{m+1}$. Then,
 ${\bbb x}^T A_0 \bbb{x}=\big(\sum_{k=0}^m x_k/2^{k}\big)^2$, and since
 $m\geq 1$, the matrix $A_0$ is singular (hence, not positive definite).
 Similarly,
 for ${\bbb x}^T=(x_0,\ldots,x_{m-1})\in \R^{m}$,
 ${\bbb x}^T B_0 \bbb{x}=(1/4)\big(\sum_{k=0}^{m-1} x_k/2^{k}\big)^2$, which
 is positive definite if and only if $m=1$ ($n=4$). On the other hand,
 $A_0(\varepsilon)=A_0\succeq0$ for all $\varepsilon\in(0,1/2)$, while
 ${\bbb x}^T B_0(\varepsilon) \bbb{x}=\big(1/4-\varepsilon(1-\varepsilon)\big)
 \big(\sum_{k=0}^{m-1} x_k/2^{k}\big)^2\geq 0$ for small enough $\varepsilon>0$.
 According to characterization (i), the given $\beta$'s are EOS, although
 the numbers $\nu_k({\bbb \beta})$ ($k=0,\ldots,n-2$)
 do not satisfy the condition $A_0\succ 0$ and $B_0\succ 0$.

 Next, suppose that $n=2m+3$. Then, $A_1=B_1$ and
 it follows that ${\bbb x}^T A_1 \bbb{x}={\bbb x}^T B_1 \bbb{x}=(1/2)
 \big(\sum_{k=0}^{m} x_k/2^{k}\big)^2$, showing that
 $A_1$ (and $B_1$) is singular and positive semi-definite.
 On the other hand,
 ${\bbb x}^T A_1(\varepsilon) \bbb{x}={\bbb x}^T B_1(\varepsilon) \bbb{x}
 =(1/2-\varepsilon)
 \big(\sum_{k=0}^{m} x_k/2^{k}\big)^2\geq$ 0, and (ii) shows that the
 $\beta$'s are EOS. Hence, although
 the numbers $\nu_k({\bbb \beta})$ ($k=0,\ldots,n-2$)
 do not satisfy the condition $A_1\succ 0$ and $B_1\succ 0$, the
 corresponding $\beta_j$ are EOS.

 It can be checked that the given $\beta$'s
 are the EOS from the two-valued r.v.\ $X$ with $\Pr(X=0)=\Pr(X=2^n)=1/2$.
 \smallskip
 \end{pr}

 \begin{REM}
 \label{rem.uniqueness}
 (a) Assume that for $i=0$ or $1$, $A_i\succeq 0$, $B_i\succeq 0$, and either
 $\det A_i=0$ or $\det B_i=0$ (or both). Then, the
 measure $\mu=\mu_0$ is $[0,1]$-determinate from its moments $(\nu_k)_{k=0}^{n-2}$;
 see Theorem 10.7 in Schm\"udgen (2017). Hence, if this is the case,
 we can find
 $\varepsilon\in(0,1/2)$ such that $A_i(\varepsilon)\succeq 0$
 and $B_i(\varepsilon)\succeq 0$ if and only if the support of
 (the unique) $\mu_0$
 does not contain any of the
 endpoints
 $0$ and $1$.
 \smallskip

 \noindent
 (b) The finite supporting set of the discrete measure $\mu_0$,
 constructed in the proof
 of Theorem \ref{theo.2},
 can be chosen to contain at most
 $k\leq n/2$ (rather than $k\leq n-1$) points.
 \end{REM}

 \begin{EXAM}
 \label{exam.n=4}
 (The case $n=4$).
 Assume we are given $\beta_1<\beta_2<\beta_3<\beta_4$.
 Since for $n=4$ we have $m=1$, the matrices
 $A_0$, $B_0$ are
 given by (see (\ref{nu_k}), (\ref{Hankel.even})),
 \[
 A_0
 \hspace*{-.4ex}
 = \frac{1}{\lambda ^2}
 \hspace*{-.4ex}
 \left(
 \hspace*{-1.3ex}
 \begin{array}{cc}
 (\beta_3-\beta_2)+3(\beta_4-\beta_1) & (\beta_4-\beta_3)+2(\beta_4-\beta_2)
 \\
 (\beta_4-\beta_3)+2(\beta_4-\beta_2) & 3 (\beta_4-\beta_3)
 \end{array}
 \hspace*{-1.3ex}
 \right),
 \ \
 B_0
 \hspace*{-.4ex}
 =
 \hspace*{-.4ex}
 \left(
 \hspace*{-1.3ex}
 \begin{array}{c}
 \frac{2}{\lambda}(\beta_3-\beta_2)
 \end{array}
 \hspace*{-1.3ex}
 \right),
  \]
with $\lambda$ as in Definition~\ref{defi.nu_k}.
  Hence, $B_0$ is (trivially) positive definite, and
  by Sylvester's criterion,
  $A_0=A_0(\varepsilon)$ is positive semi-definite if and only if
 \be
 \label{NS}
  (\beta_{2}-\beta_{1})(\beta_{4}-\beta_{3}) \geq
  \Big(\frac{2}{3} (\beta_{3}-\beta_{2})\Big)^2.
 \ee
 According to Theorem \ref{theo.2}(i), the $\beta$'s are EOS if and only if
 (\ref{NS}) is satisfied.
 Based on (\ref{NS}) we immediately deduce that, e.g.,
 the numbers $(0,\  2,\  5,\  7)$  are EOS, while the numbers $(0,\ 2,\ 11,\ 13)$ are not.
 For the first set of numbers, the r.v.\ $T$
 in (\ref{betaj binomial}) or (\ref{simple.moments}) is uniquely
 determined (in fact, $T\equiv 1/2$),
 because $\bbb{\nu}({\bbb \beta)}=(1,1/2,1/4)=(1,\E T,\E T^2)$, showing
 that $\Var T=0$; see Remark \ref{rem.uniqueness}(a) and (\ref{nu_k}) of
 Definition \ref{defi.nu_k}.
 Consequently, from Lemma \ref{lem.1} we conclude that the corresponding
 r.v.\ $X$, assuming the given expected order statistics, is also unique,
 namely, $\Pr(X=-1/2)=\Pr(X=15/2)=1/2$; see
 Remarks \ref{rem.3}, \ref{rem.2}.
 \end{EXAM}

 \begin{EXAM}
 \label{exam.n=5}
 (The case $n=5$).
 It can be checked that for $n=5$, the $2\times 2$ matrices
 $A_1$, $B_1$,
 see (\ref{Hankel.odd}), (\ref{nu_k}),
 are
 positive semi-definite if and only if
 $(\beta_3-\beta_2)(\beta_5-\beta_4)
 \geq \frac{1}{2}(\beta_4-\beta_3)^2$
 and
 $(\beta_2-\beta_1)(\beta_4-\beta_3)
 \geq \frac{1}{2}(\beta_3-\beta_2)^2$.
 Moreover, if both inequalities are strict (case $(+,+)$) then $A_1\succ 0$,
 $B_1\succ 0$, and Theorem \ref{theo.2}(iii) shows that
 the $\beta$'s are EOS. If, however, one (or both) of the inequalities
 reduces to an equality, one has to check the condition
 (ii) of Theorem \ref{theo.2} in detail. For instance,
 if both matrices are singular (case $(0,0)$),
 then $A_1(\varepsilon)\succeq 0$
 and $B_1(\varepsilon)\succeq 0$ for
 $0<\varepsilon<\min\{\beta_3-\beta_2,\beta_4-\beta_3\}/(\beta_4-\beta_2)$, and
 the $\beta$'s are again EOS. As an example of the $(0,0)$-case
 consider the numbers
 $(0,\ 1,\ 5,\ 13,\ 21)$, representing the EOS from
 the (uniquely defined) r.v.\ $X$
 with $\Pr(X=-1/10)=2/3$, $\Pr(X=121/5)=1/3$. However, both cases
 $(0,+)$ (e.g., $(0,\ 9,\ 11,\ 13,\ 14)$) and $(+,0)$
 (e.g., $(0,\ 1,\ 3,\ 5,\ 14)$)
 imply that the $\beta$'s are not EOS.
 To see this,
 assume that $2 (\beta_3-\beta_2)(\beta_5-\beta_4)
 = (\beta_4-\beta_3)^2$ and
 $2(\beta_2-\beta_1)(\beta_4-\beta_3)>(\beta_3-\beta_2)^2$; case $(0,+)$. Then,
 $B_1\succ 0$ (hence, $B_1(\varepsilon)\succeq 0$ for small $\varepsilon>0$) and
 $A_1\succeq 0$ with $\det A_1=0$. It can be verified that for
 ${\bbb{x}}^T=(x_0,x_1):=(\beta_4-\beta_3,-(\beta_4-\beta_2))$,
 $\bbb{x}^T A_1(\varepsilon){\bbb x}=-\varepsilon \Delta$, where
 $\Delta>0$
 depends only on
 $\beta$'s, and thus, according to Theorem \ref{theo.2}(ii),
 the $\beta$'s cannot be EOS. By the same reasoning, this is also true
 for the $(+,0)$-case. Therefore, the complete characterization for $n=5$
 says for the $\beta$'s to be EOS it is necessary and sufficient that
 either
 $2 (\beta_2-\beta_1)(\beta_4-\beta_3)
 = (\beta_3-\beta_2)^2$ and
 $2 (\beta_3-\beta_2)(\beta_5-\beta_4)= (\beta_4-\beta_3)^2$,
 or
 $2 (\beta_2-\beta_1)(\beta_4-\beta_3)
 > (\beta_3-\beta_2)^2$ and
 $2 (\beta_3-\beta_2)(\beta_5-\beta_4)> (\beta_4-\beta_3)^2$;
 that is,
 either both matrices
 $A_1$, $B_1$ are positive definite, or both are positive semi-definite
 and singular.
 We do not know if the situation is similar for
 odd values of $n\geq 7$.
 \end{EXAM}

 Our final result characterizes the binomial mixtures for which
 the mixing distribution is supported in the open interval
 $(0,1)$ (cf.\ Wood, 1992, 1999). The proof,
 being an immediate application of Theorems \ref{theo.1},
 \ref{theo.2}, is omitted.

 \begin{theo}
 \label{theo.convexhul}
 Let ${\bbb p}=(p_0,\ldots,p_n)$ {\rm ($n\geq 2$)} be a probability
 vector
 {\rm ($p_i\geq 0$,
 $\sum_{i=0}^n p_i=1$)} and ${\bbb u}={\bbb u}({\bbb p})
 =(u_0,\ldots,u_n)$, where
 \[
 u_k:=
 {{n\choose k}}^{-1} \sum_{j=k}^n {j\choose k} p_j,
 \ \ k=0,1,\ldots,n.
 \]
 If $n=2m$ set
 \[
 A(\varepsilon)=\Big(u_{i+j}\Big)_{i,j=0}^m, \ \ \
 B(\varepsilon)=\Big(u_{i+j+1}-u_{i+j+2}-
 \varepsilon(1-\varepsilon)u_{i+j}\Big)_{i,j=0}^{m-1},
 \]
 and if $n=2m+1$ set
 \[
 A(\varepsilon)=\Big(u_{i+j+1}-\varepsilon u_{i+j}\Big)_{i,j=0}^m, \ \ \
 B(\varepsilon)=\Big((1-\varepsilon)u_{i+j}-u_{i+j+1}\Big)_{i,j=0}^{m}.
 \]
 Then, the following
 are equivalent.
 \smallskip

 \noindent {\rm (i)} $A(\varepsilon)\succeq 0$ and $B(\varepsilon)\succeq 0$
 for some $\varepsilon$ with $0<\varepsilon<1/2$.
 \smallskip

 \noindent {\rm (ii)} ${\bbb p}\in \mbox{\rm Conv}[B_0]$, where
 \[
 B_0=\left\{\left({n\choose j}p^j(1-p)^{n-j}\right)_{j=0}^n, \ 0<p<1 \right\}
 \]
 is the
 open
 binomial probability curve {\rm (without its endpoints)} and
 $\mbox{\rm Conv}[X]$ denotes the convex hull of $X\subseteq \R^{n+1}$.
 \smallskip

 \noindent {\rm (iii)} There exists an r.v.\ $V$ with $\Pr(0<V<1)=1$
 such that
 \[
 p_j=\E\left\{ {n\choose j} V^j (1-V)^{n-j}\right\}, \ \ \ j=0,1,\ldots,n.
 \]
 \smallskip

 \noindent {\rm (iv)} ${\bbb u}\in \mbox{\rm Conv}[M_0]$, where
 $M_0=\left\{(1,t,t^2,\ldots,t^n), \ \ 0<t<1 \right\}$
 is the
 open
 moment curve {\rm (without its endpoints)}.
 \smallskip

 \noindent {\rm (v)} There exists an r.v.\ $V$ with $\Pr(0<V<1)=1$
 such that
 \[
 \medskip
 u_k=\E V^k, \ \ \ k=0,1,\ldots,n.
 \]
 \end{theo}

 Let ${\bbb{x}}(t)=(1,t,t^2,\ldots,t^n)$, $0\leq t\leq 1$.
 A simple application of Theorem \ref{theo.convexhul} shows
 that for $n\geq 3$ (in contrast to the case $n=2$),
 the line segment $(1-\lambda){\bbb{x}(0)}+\lambda{\bbb{x}(t_0)}$,
 $0\leq \lambda<1$, lies outside $\mbox{\rm Conv}[M_0]$.
 Given $\lambda_0,t_0\in (0,1)$, it follows from Farkas' Lemma
 (see, e.g., Bertsimas and Tsitsiklis (1997), Theorem 4.6)
 that for any $m$,
 and any given collection
 $\{t_0,\ldots,t_m\}\subseteq (0,1)$, we can find a polynomial
 $p$ with $\deg(p)\leq n$, such that
 $(1-\lambda_0)p(0)+\lambda_0 p(t_0)<0$
 and $p(t_i)\geq 0$,
 $i=0,1,\ldots,m$.

 {

 }

 {

 \footnotesize

 \vspace*{2em}

 \appendix
 \section{}
 \label{sec.app.a}
 \begin{center}
 \bf
   A detailed proof of Lemma \ref{lem.1}
 \end{center}

 \noindent
 Provided that (\ref{EXkk}) holds true, the uniqueness of $X$ follows from
 the classical result of Hoeffding (1953), since the sequence $\E X_{k:k}$
 determines the triangular array $\E X_{j:n}$, and vice-versa;
 see (\ref{mujn from muk}), (\ref{muk from mjn}).
 We proceed to verify that $F_0^{-1}$ is indeed the distribution
 inverse of an integrable r.v.\ $X$ satisfying (\ref{EXkk}).
 Left continuity of $F_0^{-1}$ follows automatically from
 its definition. Moreover,
 the function $g(t):=F_0^{-1}(t)+c_T$ can be written as
 \be
 \label{W}
 \hspace*{-.3ex}
 g(t)
 \hspace*{-.5ex}=
 \hspace*{-.7ex}
 \left\{
 \hspace*{-1.5ex}
 \begin{array}{cc}
 -4\big[F_T(\mbox{$\frac{1}{2}-$})-F_T(t-)\big]
 -\int_{t}^{1/2} \frac{1-2u}{u^2(1-u)^2} \big[F_T(u)-F_T(t-)\big] d u,
 &
 \hspace*{-.5ex}
 0<t\leq \frac{1}{2},
 \\
 4\big[F_T(t-)-F_T(\mbox{$\frac{1}{2}-$})\big]
 +\int_{1/2}^{t} \frac{2u-1}{u^2(1-u)^2} \big[F_T(t-)-F_T(u)\big] d u,
 &
 \hspace*{-.5ex}
 \frac{1}{2}\leq t<1,
 \end{array}
 \right.
 \ee
 showing that $g(t)\geq 0$ for $t\geq 1/2$ and $g(t)\leq 0$ for $t\leq 1/2$.
 For $t_1<t_2$ with $1/2\leq t_1<t_2<1$ we have
 \[
 g(t_2)-g(t_1)=
 \frac{F_T(t_2-)-F_T(t_1-)}{t_1(1-t_1)}
 +
 \int_{t_1}^{t_2}
 \frac{2u-1}{u^2(1-u)^2}\big[F_T(t_2-)-F_T(u)\big] du\geq 0.
 \]
 Similarly, for $t_1<t_2$ with $0<t_1<t_2\leq 1/2$,
 \[
 g(t_2)-g(t_1)
 =
 \frac{F_T(t_2-)-F_T(t_1-)}{t_2(1-t_2)}
 +
 \int_{t_1}^{t_2}
 \frac{1-2u}{u^2(1-u)^2}\big[F_T(u)-F_T(t_1-)\big] du\geq 0.
 \]
 Therefore, $g$, and hence $F_0^{-1}$, is nondecrasing. In order to verify that
 $F_0^{-1}\in L^1$, we shall calculate the integrals
 \[
 J_1=\int_0^{1/2} |g(t)| d t= \int_0^{1/2} -g(t) d t, \ \
 J_2=\int_{1/2}^1 |g(t)| d t  = \int_{1/2}^1 g(t) d t.
 \]
  In the following we shall make repeatedly
  use of the fact that for a nonegative
  r.v.\ $Y$, $\E Y=\int_{0}^{\infty}\Pr(Y>t)dt$
  or $\E Y=\int_{0}^{\infty}\Pr(Y\geq t)dt$, followed
  by subsequent applications of Tonelli's theorem.
 Using (\ref{W}) and Tonelli's theorem we have
 \[
 J_1=2F_T(\mbox{$\frac{1}{2}-$})-4\int_{0}^{1/2} F_T(t) dt
 +\int_0^{1/2} \frac{1-2u}{u^2(1-u)^2} \int_0^u \big[ F_T(u)-F_T(t)\big] dt du,
 \]
 on noting that $F_T(t)=F_T(t-)$ a.e. Considering the nonnegative r.v.\
 $Y=h(T)=T I(T\leq u)$
 it is easily seen that $\Pr(Y>t)=F_T(u)-F_T(t)$ for $t<u$, and the probability is zero for
 $t\geq u$.
 Hence, $\E Y=\int_0^u \big[ F_T(u)-F_T(t)\big] dt$, and also,
 $\E h(T)= \int_{(0,u]} t
 dF_T(t)$.
 Since these expectations are equal, we obtain
 \[
 \int_0^u \big[ F_T(u)-F_T(t)\big] dt = \int_{(0,u]} t dF_T(t).
 \]
 Substituting this equality to the double integral in $J_1$ and interchanging once again
 the order of
 integration (since the integrand is nonnegative), we obtain
 \[
 \int_0^{1/2} \frac{1-2u}{u^2(1-u)^2} \int_0^u \big[ F_T(u)-F_T(t)\big] dt du
 =\int_{(0,1/2]} \big(\frac{1}{1-t}-4t\big) dF_T(t)<\infty.
 \]
 In order to evaluate the exact value of $J_1$, it remains to express the
 integral $\int_{0}^{1/2} F_T(t) dt$ in terms of integrals w.r.t.\
 $dF_T$.
 For $u\in(0,1)$,
 consider the nonnegative r.v.\ $Y=h(T)=T I(T<u)$, for which $\Pr(Y>t)=F_{T}(u-)-F_{T}(t)$
 for $t<u$, and zero otherwise. Then, $\E Y=\int_{0}^{u} \big[F_{T}(u-)-F_{T}(t)\big] dt
 =u F_T(u-)-\int_{0}^u F_T(t) dt$, and $\E h(T)=\int_{(0,u)}t d F_T(t)$; thus,
 \[
 \int_{0}^{u} F_T(t) dt = u F_T(u-)-\int_{0}^{u} \big[F_T(u-)-F_T(t)\big] dt
 = u F_T(u-)-\int_{(0,u)}t d F_T(t).
 \]
 Setting $u=1/2$ we find
 \[
 2F_T(\mbox{$\frac{1}{2}-$})-4\int_{0}^{1/2} F_T(t) dt= \int_{(0,1/2)} 4t d F_T(t),
 \]
 and finally, since $\int_{(0,1/2]} \big(\frac{1}{1-t}-4t\big) dF_T(t)=
 \int_{(0,1/2)} \big(\frac{1}{1-t}-4t\big) dF_T(t)$ (because the integrand vanish for $t=1/2$),
 we conclude that
 \[
 J_1=\int_{(0,1/2)} \frac{1}{1-t} d F_T(t)=\E\Bigg[\frac{1}{1-T}I(T<1/2)\Bigg].
 \]
 Using (\ref{W}) we rewrite $J_2$ as
 \[
 J_2=-2F_T(\mbox{$\frac{1}{2}-$})+4\int_{1/2}^{1} F_T(t) dt
 +\int_{1/2}^1 \frac{2u-1}{u^2(1-u)^2} \int_u^1 \big[ F_T(t)-F_T(u)\big] dt du,
 \]
 by Tonelli's theorem. Substituting
 \[
 \int_u^1 \big[ F_T(t)-F_T(u)\big] dt=\int_{(u,1)} (1-t) d F_T(t)
 \]
 in the inner integral (noting that both integrals represent the expectation
 of $Y=(1-T)I(T>u)$), and changing the order of integration,  we arrive at
 \[
 J_2=-2F_T(\mbox{$\frac{1}{2}-$})+4\int_{1/2}^{1} F_T(t) dt
 +\int_{(1/2,1)} \Big(\frac{1}{t}-4(1-t)\Big) d F_T(t)<\infty.
 \]
 For $u\in(0,1)$ consider the r.v.\ $Y=(1-T)I(T\geq u)$,
 so that, $\Pr(Y\geq y)=F_T(1-y)-F_T(u-)$
 for $y\leq 1-u$, and zero otherwise. Then,
 \[
 \E Y =\int_{0}^{1-u} \big[F_T(1-y)-F_T(u-)\big] d y
 = \int_{u}^{1} \big[F_T(t)-F_T(u-)\big] d t,
 \]
 and this expectation is also equal to $\int_{[u,1)}(1-t) d F_T(t)$. Hence,
 \[
 \int_{u}^{1} F_T(t) d t=(1-u) F_T(u-)+\int_{[u,1)}(1-t) d F_T(t).
 \]
 Substituting $u=1/2$ we obtain
 \[
 -2F_T(\mbox{$\frac{1}{2}-$})+4\int_{1/2}^{1} F_T(t) dt=
 \int_{[1/2,1)}4(1-t) d F_T(t),
 \]
 and since $\int_{(1/2,1)} \big(\frac{1}{t}-4(1-t)\big) d F_T(t)
 =\int_{[1/2,1)} \big(\frac{1}{t}-4(1-t)\big) d F_T(t)$,
 we conclude that
 \[
 J_2=\int_{[1/2,1)} \frac{1}{t} d F_T(t)=\E\Bigg[\frac{1}{T}I(T\geq 1/2)\Bigg].
 \]
 The preceding argument not only shows that $F_0^{-1}\in L^{1}$, but also proves that
 \[
 \int_{0}^1 g(t) dt=J_2-J_1=c_T,
 \]
 with $c_T$ as in (\ref{cT}), and therefore, the expectation of the
 r.v.\ $X$ with inverse d.f.\ $F_0^{-1}=g-c_T$ is zero.

 Next, set $\mu_k=E X_{k:k}=k\int_0^1 t^{k-1}F_0^{-1}(t)dt$, $k=1,2,\ldots$, and
 let $R=F_T(\mbox{$\frac{1}{2}-$})$.
 In view of  (\ref{W}), write
 $\mu_k+c_T=\int_0^1 k t^{k-1} g(t) dt=I_2-I_1$ where
 \begin{eqnarray*}
 I_1
 \hspace*{-1.5ex}&=&\hspace*{-1.5ex}
 \int_{0}^{1/2} 4k t^{k-1}
 \big[
 R
 -F_T(t)\big] dt + \int_0^{1/2}
 \int_{t}^{1/2} \frac{(1-2u)kt^{k-1}}{u^2(1-u)^2} \big[F_T(u)-F_T(t)\big] du dt,
 \\
 I_2
 \hspace*{-1.5ex}&=&\hspace*{-1.5ex}
 \int_{1/2}^{1} 4kt^{k-1}
 \big[F_T(t)-
 R
 \big] dt + \int_{1/2}^1
 \int_{1/2}^t \frac{(2u-1)kt^{k-1}}{u^2(1-u)^2} \big[F_T(t)-F_T(u)\big] du dt,
 \end{eqnarray*}
 noting that the integrands may be different form the original ones (suggested from
 (\ref{W})) at sets of measure zero.
 Changing the order of integration, according to Tonelli's theorem, we see that
 \begin{eqnarray*}
 I_1
 \hspace*{-1.5ex}&=&\hspace*{-1.5ex}
 \int_{0}^{1/2} 4kt^{k-1}\big[
 R
 -F_T(t)\big] dt + \int_0^{1/2}
 \frac{1-2u}{u^2(1-u)^2}
 \int_{0}^{u} kt^{k-1}  \big[F_T(u)-F_T(t)\big] dt du,
 \\
 I_2
 \hspace*{-1.5ex}&=&\hspace*{-1.5ex}
 \int_{1/2}^{1} 4k t^{k-1}
 \big[F_T(t)-
 R
 \big] dt
 + \int_{1/2}^1 \frac{2u-1}{u^2(1-u)^2}
 \int_{u}^1 kt^{k-1} \big[F_T(t)-F_T(u)\big] dt du.
 \end{eqnarray*}
 Consider the expectation of $Y=T^k I(T\leq u)$. Since $\Pr(Y>y)=F_T(u)-F_T(y^{1/k})$
 for $y<u^k$, and zero otherwise, we obtain
 \[
  \int_{(0,u]} t^k d F_T(t)=\E Y = \int_0^{u^k} \big[F_T(u)-F_T(y^{1/k})\big] d y
 = \int_0^{u} k t^{k-1}\big[F_T(u)-F_T(t)\big] d t.
 \]
 Thus,
 \begin{eqnarray*}
 I_1
 \hspace*{-1.3ex}&=&\hspace*{-1.3ex}
 \int_{0}^{1/2} 4k t^{k-1}\big[F_T(\mbox{$\frac{1}{2}-$})-F_T(t)\big] dt
 +
 \int_0^{1/2}
 \frac{1-2u}{u^2(1-u)^2}
 \int_{(0,u]} t^k d F_T(t) du
 \\
 \hspace*{-1.3ex}&=&\hspace*{-1.3ex}
 \int_{0}^{1/2} 4k t^{k-1}\big[F_T(\mbox{$\frac{1}{2}-$})-F_T(t)\big] dt
 +
 \int_{(0,1/2]} t^k\Big(\frac{1}{t(1-t)}-4\Big) d F_T(t).
 \end{eqnarray*}
 Next, set $Y=T^k I(T<u)$ (for $0<u<1$) with $\E Y= \int_{(0,u)} t^k d F_T(t)$,
 and
 observe
 that $\Pr(Y>y)=F_T(u-)-F_T(y^{1/k})$ for $y<u^k$ (and zero otherwise), to conclude
 the identity
 \[
 \int_{(0,u)} t^k d F_T(t)=\int_0^{u^k} \big[F_T(u-)-F_T(y^{1/k})\big] dy
 =
 \int_0^{u} k t^{k-1}\big[F_T(u-)-F_T(t)\big] dt.
 \]
 Applying this with $u=1/2$ we obtain an explicit simple formula for $I_1$:
 \[
 I_1=\int_{(0,1/2)} 4t^k d F_T(t)
 +
 \int_{(0,1/2]} t^k\Big(\frac{1}{t(1-t)}-4\Big) d F_T(t)=
 \int_{(0,1/2)} \frac{t^k}{t(1-t)} d F_T(t).
 \]

 Finally, in order to calculate $I_2$, consider the auxiliary variable $Y=(1-T^k)I(T>u)$
 (with $0<u<1$), for which $\Pr(Y\geq y)=F_T((1-y)^{1/k})-F_T(u)$ for $y<1-u^k$, and
 zero otherwise. The alternative expressions for its expectation yield
 \[
 \int_{(u,1)} \hspace*{-1.4ex}(1-t^k) d F_T(t)= \int_0^{1-u^k}
 \hspace*{-1.2ex}
 \big[F_T((1-y)^{1/k})-F_T(u)\big] dy =
 \int_u^{1} \hspace*{-1.2ex} k t^{k-1}\big[F_T(t)-F_T(u)\big] dt.
 \]
 Hence,
 \[
 I_2
 =
 \int_{1/2}^{1} 4k t^{k-1} \big[F_T(t)-F_T(\mbox{$\frac{1}{2}-$})\big] dt
 + \int_{1/2}^1 \frac{2u-1}{u^2(1-u)^2}
 \int_{(u,1)} (1-t^k) d F_T(t) du,
 \]
 and applying once again Tonelli's theorem, we get
 \[
 I_2
 =
 \int_{1/2}^{1} 4k t^{k-1} \big[F_T(t)-F_T(\mbox{$\frac{1}{2}-$})\big] dt
 +
 \int_{(1/2,1)} (1-t^k)\Big(\frac{1}{t(1-t)}-4\Big) d F_T(t).
 \]
 If we set $Y=(1-T^k)I(T\geq u)$, we see that $\Pr(Y\geq y)=F_T((1-y)^{1/k})-F_T(u-)$
 for $y\leq 1-u^k$, and zero otherwise, obtaining
 \[
 \int_{[u,1)}
 \hspace*{-1.3ex}
 (1-t^k) d F_T(t)
 \hspace*{-.4ex}
 =
 \hspace*{-.7ex}
 \int_0^{1-u^k}
 \hspace*{-1ex}
 \big[F_T((1-y)^{1/k})-F_T(u-)\big] dy
 \hspace*{-.4ex}
 =
 \hspace*{-.7ex}
 \int_u^{1}
 \hspace*{-1ex}
 k t^{k-1} \big[F_T(t)-F_T(u-)\big]
 \hspace*{-.3ex}
 dt.
 \]
 Applying this identity with $u=1/2$ we conclude that
 \[
 I_2
 =
 \hspace*{-.4ex}
 \int_{[1/2,1)}
 \hspace*{-2ex}
 4(1-t^k) d F_T(t) +
 \int_{(1/2,1)}
 \hspace*{-2ex}
 (1-t^k)\Big(\frac{1}{t(1-t)}-4\Big) d F_T(t)
 =
 \hspace*{-.4ex}
 \int_{[1/2,1)}
 \hspace*{-.5ex}
 \frac{1-t^k}{t(1-t)} d F_T(t).
 \]
 By the preceding calculations,
 \[
 \mu_k+c_T=I_2-I_1=\int_{[1/2,1)} \frac{1-t^k}{t(1-t)} d F_T(t
 )-
 \int_{(0,1/2)} \frac{t^k}{t(1-t)} d F_T(t)
 \]
 (observe that the r.h.s.\ equals to $c_T$ for $k=1$, showing once again that $\mu_1=0$).
 Therefore, for $k=0,1,\ldots$,
 \begin{eqnarray*}
 \mu_{k+2}-\mu_{k+1}
 \hspace*{-1.5ex}&=&\hspace*{-1.5ex}
 \int_{[1/2,1)} \frac{(1-t^{k+2})-(1-t^{k+1})}{t(1-t)} d F_T(t)-
 \int_{(0,1/2)} \frac{t^{k+2}-t^{k+1}}{t(1-t)} d F_T(t)
 \\
 \hspace*{-1.5ex}&=&\hspace*{-1.5ex}
 \int_{(0,1)} t^k d F_T(t)
 \\
 \hspace*{-1.5ex}&=&\hspace*{-1.5ex}
 \E T^k, \Bigg.
 \vspace*{-1em}
 \end{eqnarray*}
 and the proof is complete.
 \vspace*{1.2em}

 }


\begin{thebibliography}{99}

 \bibitem{ABN1992}
 {\sc Arnold, B.C.; Balakrishnan, N.; Nagaraja, H.N.}\ (1992).
 \textit{A First Course in
 Order Statistics.}
 John Wiley \& Sons, New York.

 \bibitem{B-T1997}
 {\sc Bertsimas, D.; Tsitsiklis, J.N.}\ (1997).
 \textit{Introduction to Linear Optimization.}
 Athena Scientiffic, Belmont, Massachusetts.

 \bibitem{DN2003}
 {\sc David, H.A.; Nagaraja, H.N.}\ (2003).
 \textit{Order Statistics} (3rd ed.),
 John Wiley \& Sons, Hoboken, New Jersey.

 \bibitem
 {Hausdorff1921}
    {\sc Hausdorff, F.}\ (1921).
    Summationmethoden und momentfolgen. I.
    \textit{Math.\ Zeitchrift},  {\bf9}(1), 74--109.

 \bibitem
 {Hoeffding1953}
    {\sc Hoeffding, W.}\ (1953).
    On the distribution of the expected values of the order statistics.
    \textit{Ann.\ Math.\ Statist.},  {\bf24}(1), 93--100.

 \bibitem
 {Huang1998}
    {\sc Huang, J.S.}\ (1998).
    Sequences of expectations
    of maximum-order statistics.
    \textit{Statist.\ Probab.\ Lett.}, {\bf38}, 117--123.

 \bibitem{J+B}
 {\sc Jones, M.C.; Balakrishnan, N.}\ (2002).
 How are moments and moments
 of spacings related to distribution functions?
 {\it J.\ Stat.\ Plann.\ Inference}
 (C.R.\ Rao 80th birthday felicitation volume,
 Part I), {\bf 103},
  377--390.

 \bibitem
 {Kadane1971}
    {\sc Kadane, J.B.}\ (1971).
    A moment problem for order statistics.
    \textit{Ann.\ Math.\ Statist.},  {\bf42},
    745--751.

 \bibitem
 {Kadane1974}
    {\sc Kadane, J.B.}\ (1974).
    A characterization of triangular arrays which are
    expectations of order statistics.
    \textit{J.\ Appl.\ Probab.},  {\bf11},
    413--416.

 \bibitem
 {KS1966}
    {\sc Karlin, S.; Studden, W.}\ (1966). \textit{Tchebycheff Systems
    With Applications
    in Analysis and Statistics}. Interscience, New York.

 \bibitem
 {Kolo2000}
    {\sc Kolodynski, S.}\ (2000).
    A note on the sequence of expected extremes.
    \textit{Statist.\ Probab.\ Lett.}, {\bf47}, 295--300.

 \bibitem
 {Mallows1973}
    {\sc Mallows, C.L.}\ (1973).
    Bounds on distribution functions in terms of expectations
    of order-statistics.
    \textit{Ann.\  Probab.}, {\bf1}, 297--303.

 \bibitem
  {Pap2001}
  {\sc Papadatos, N.}\ (2001).
  Distribution and expectation bounds on order statistics from possibly dependent variates.
  \textit{Statist.\ Probab.\ Lett.},  {\bf 54},
  21--31.

 \bibitem
  {Pap2017}
   {\sc Papadatos, N.}\ (2017).
   On sequences of expected maxima and expected ranges.
   \textit{J.\ Appl.\ Probab.}, {\bf 54}, 1144--1166.

 \bibitem
  {Paarson1902}
   {\sc Pearson, K.}\ (1902).
   Note on Francis Galton's problem.
   \textit{Biometrika}, {\bf 1}, 390--399.

 \bibitem
 {Schmudgen2017}
 {\sc Schm\"{u}dgen, K.}\ (2017). {\it The Moment Problem}.
 Graduate Texts in Mathematics, {\bf 277}, Springer.

 \bibitem
 {Wood1992}
 {\sc Wood, G.R.}\ (1992). Binomial mixtures and finite exchangeability.
   \textit{Ann.\ Probab.},
   {\bf 20}(3), 1167--1173.

 \bibitem
 {Wood1999}
 {\sc Wood, G.R.}\ (1999). Binomial mixtures: geometric estimation
  of the mixing distribution.
   \textit{Ann.\ Statist.},
   {\bf 27}(5), 1706--1721.
 \end{thebibliography}
 \end{document}